\documentclass[12pt]{amsart}
\usepackage[T1]{fontenc}
\usepackage{amssymb}
\let\SavedRightarrow=\Rightarrow
\usepackage{marvosym}
\let\Rightarrow=\SavedRightarrow

\newcommand{\Aaa }{\mathcal A}
\newcommand{\Raa }{\mathcal R}
\newcommand{\Bee }{\mathcal B}
\newcommand{\Cee }{\mathcal C}

\newcommand{\Qee }{\mathcal Q}
\newcommand{\Pee }{\mathcal P}

\newcommand{\Tee }{\mathcal T}
\newcommand{\See }{\mathcal S}

\newcommand{\cl}{\operatorname{cl}}
\renewcommand{\int}{\operatorname{Int}}

\newtheorem{theorem}{Theorem}
\newtheorem{corollary}[theorem]{Corollary}
\newtheorem{lemma}[theorem]{Lemma}

\newtheorem{proposition}[theorem]{Proposition}

\author{Andrzej Kucharski}
\address{Andrzej Kucharski \\
 Institute of Mathematics, University of
Silesia \\
 ul. Bankowa 14, 40-007 Katowice}
\email{akuchar@ux2.math.us.edu.pl}
\author{Szymon Plewik}
\address{Szymon Plewik\\Institute of Mathematics,
University of Silesia, ul. Ban\-ko\-wa 14, 40-007 Katowice}
\email{plewik@ux2.math.us.edu.pl}

\begin{document}

\title{Inverse systems and I-favorable spaces} 

\begin{abstract}
We show that a  compact  space  is I-favorable if, and only if it can be represented as the limit of a $\sigma$-complete inverse system  of  compact metrizable spaces  with skeletal bonding maps. We also show that any completely regular I-favorable space can be embedded as a dense subset of the limit of a $\sigma$-complete inverse system  of  separable metrizable spaces  with skeletal bonding maps. 
\end{abstract}

\maketitle
\section{Introduction} 
We investigate  the class of all limits
 of $\sigma$-complete inverse systems  of  compact  metrizable spaces  with skeletal bonding maps. Notations are used  the same as in the monograph \cite{eng}. For example,  a compact space is Hausdorff, and a regular space is $T_1$. A directed set $\Sigma $ is said to be
$\sigma$-complete if any countable chain of its elements has
least upper bound in $\Sigma$. 
 An inverse system $  \{ X_\sigma , \pi^\sigma_\varrho, \Sigma\}$ is said to be a $\sigma$-complete, whenever  $\Sigma $ is 
$\sigma$-complete and for every chain $\{\sigma_n: n \in \omega \} \subseteq \Sigma$, such that $\sigma = \sup \{\sigma_n: n \in \omega \} \in \Sigma,$ there  holds $$X_{\sigma }= \varprojlim \{ X_{\sigma_n}, \pi^{\sigma_{n+1}}_{\sigma_n}\},$$ compare  \cite{s81}. However, we will consider inverse systems where bonding maps are surjections. Another  details about inverse systems one can find in \cite{eng} pages 135 - 144. For  basic facts about I-favorable spaces we refer to \cite{dkz}, compare also \cite{kp}.

Through the course of this note we  modify quotient topologies and quotient maps,  introducing $\Qee_{\Pee}$-topologies and $\Qee_{\Pee}$-maps, where $\Pee$ is a family of subsets of $X$.   Next,  we assign  the family   
$\Pee_{seq}$ (of all sets with  some properties of cozero sets) to a given family $\Pee$. Frink's theorem is used to show that the  $\Qee_{\Pee}$-topology is  completely regular,  whenever $\Pee \subseteq \Pee_{seq}$ is a ring of subsets of $X$, see Theorem  \ref{t5}. Afterwards,  some special club filters are described as  systems of countable skeletal families. This yields that each family which belongs to a such club filter is a countable skeletal family, which  produces a skeletal map onto a compact metrizable space. Theorem \ref{10} is the main result: I-favorable compact  spaces  coincides with limits of  $\sigma$-complete inverse systems  of  compact metrizable spaces  with skeletal bonding maps. 

 E.V. Shchepin has considered several classes of compact spaces in a few papers, for example   \cite{s76}, \cite{s79} and \cite{s81}.  He introduced the class of compact openly generated  spaces.  A compact space $X$ is called \textit{openly generated}, whenever $X$  is the limit
 of a $\sigma$-complete inverse system   of compact metrizable spaces  with open bonding maps. Originally, Shchepin used  another name: open-generated spaces; see \cite{s81}. A. V. Ivanov showed that a compact space $X$ is openly generated if, and only if its superextension is a Dugundji space, see \cite{iva}. Then 
 Shchepin established that the classes of openly generated compact
 spaces and of $\kappa$-metrizable spaces are the same, see Theorem 21 in \cite{s81}. Something likewise is    established for compact I-favorable spaces in Theorem \ref{10}.

A Boolean algebra $\mathbb B$ is semi-Cohen (regularly filtered) if, and only if $[\mathbb B]^\omega$ has a closed unbounded set of countable regular subalgebras, in other words $[\mathbb B]^\omega$ contains a club filter. Hence, the Stone space of a semi-Cohen algebras is  I-favorable. Translating Corollary 5.5.5 of L. Heindorf and L. B. Shapiro  \cite{hs}   on topological notions, one can obtain our's main result in zero-dimensional cases,  compare also Theorem 4.3 of B. Balcar, T. Jech and J. Zapletal   \cite{bjz}. We get  Theorem \ref{9} which says that each completely regular  I-favorable space is homeomorphic to a dense subspace of the limit of an inverse system $\{ X/\Raa, q^\Raa_\Pee,\Cee\}$, where spaces $X/\Raa$ are metrizable and separable, bonding maps $q^\Raa_\Pee$ are skeletal and the directed set $\Cee$ is $\sigma$-complete.

 \section{$\Qee_\Pee$-topologies}  
  Let $\Pee$ be a family  of subsets of $X$.  We say that $y\in [x]_{\Pee}$, whenever   $x\in V$ if, and only if $y\in V$, for each $V\in \Pee$. The family of all  classes $[x]_{\Pee}$  is denoted  $X{/\Pee}$. Note that   $[x]_\Pee \subseteq V$ if, and only if $[x]_\Pee \cap V\not=\emptyset$, for each $V\in \Pee$. Put $q(x)=[x]_\Pee$.  The function $q:X \to X/\Pee$ is called an \textit{$\Qee_\Pee$-map}. The coarser topology on $X{/\Pee}$ which contains all  images $q[V] = \{ [x]_\Pee: x \in V \}$, where $V\in \Pee$,  is called an \textit{$\Qee_\Pee$-topology}. 
 If $V\in \Pee$, then $q^{-1}(q[V]) = V$. 
Indeed, we have $V\subseteq q^{-1}(q[V])$, since $q: X\to X{/\Pee}$ is a surjection. Suppose $x\in q^{-1}(q[V])$. Then  $q(x)\in q[V]$, and $[x]_\Pee \cap V\not=\emptyset$.  We get $[x]_\Pee \subseteq V$, since $V\in \Pee$. Therefore   $x\in V$.

\begin{lemma}\label{tie} Let $\Pee$ be a family of open subsets of  a topological space $X$. If $\Pee$ is a closed under finite intersections, then the $\Qee_\Pee$-map $q:X\to X{/\Pee} $ is continuous. Moreover, if  $X=\bigcup\Pee$, then   the family $\{q[V]: V\in \Pee\}$ is a base for the $\Qee_\Pee$-topology.
\end{lemma}
\begin{proof} We have $q[V\cap U]  = q[V]\cap q[U]$, for every $U, V\in \Pee$. Hence, the family $\{ q[V]: V \in \Pee\}$ is closed under finite intersections. This family is  a base for the $\Qee_\Pee$-topology, since $X=\bigcup\Pee$ implies that $X{/\Pee}$ is an union of basic sets. Obviously, the $\Qee_\Pee$-map $q$ is continuous.  
\end{proof}
  
 Additionally, if $X$ is a compact space and $X/\Pee$ is  Hausdorff, then the $\Qee_\Pee$- map $q: X\to X{/\Pee}$ 
 is a quotient map. Also, the  $\Qee_\Pee$-topology coincides with the quotient topology, compare \cite{eng} p. 124.

 Let $\Raa$ be 
a family of  subsets of $X$. Denote by $\Raa_{seq}$ the family of all sets $W$ which  satisfy the  following condition:  \textit{There exist sequences  $ \{ U_n: n \in \omega\} \subseteq \Raa$ and $ \{ V_n: n \in \omega\} \subseteq \Raa$ such that $U_k\subseteq (X\setminus V_k)  \subseteq U_{k+1}$, for any $k\in\omega$, and $\bigcup \{ U_n: n \in \omega\}=W$}. 

If   $\Raa_{seq}\not= \emptyset$, then  $\bigcup\Raa=X$. Indeed, take $W\in \Raa_{seq}$. Whenever $U_n$ and $V_n$ are elements of sequences witnessing $W \in \Raa_{seq}$, then $X \setminus V_k \subseteq U_{k+1} \subseteq W$ implies $U_{k+1} \cup V_k = X$.

If $X$ is   a completely  regular space and $\Tee$ consists of all cozero sets of $X$, then $\Tee=\Tee_{seq}$. Indeed, for each $W\in \Tee$, fix a  continuous function $f: X \to [0,1]$ 
such that $W=f^{-1}((0,1])$.  Put  $U_n = f^{-1}((\frac 1 n,1])$ and $X \setminus V_n = f^{-1}([\frac 1 n,1])$. 

 Recall that, a family of sets is called a \textit{ring of sets} whenever it is closed under finite intersections and finite unions.
\begin{lemma}\label{4hau} If a ring of sets $\Raa$ is contained in  $\Raa_{seq}$, then any countable union $ \bigcup \{U_n \in \Raa: n\in \omega \}$ belongs to  $\Raa_{seq}$. 
\end{lemma}
\begin{proof} Suppose that sequences  $\{ U^n_k: k\in \omega\}\subseteq \Raa$  and $\{ V^n_k: k\in \omega\}\subseteq \Raa$  witnessing   $U_n \in \Raa_{seq}$, respectively. Then sets $U_n^0\cup U_n^1 \cup \ldots \cup U_n^n$ and $V_n^0\cap V_n^1 \cap \ldots \cap V_n^n$ are successive elements of sequences which  witnessing  $ \bigcup \{U_n \in \Raa: n\in \omega \} \in \Raa_{seq}$.
 \end{proof}
 \begin{lemma}\label{h5}  If a family of sets $\Pee$ is contained in $ \Pee_{seq}$, then the  \mbox{ $\Qee_\Pee$-topology}  is   Hausdorff.
\end{lemma}
\begin{proof} Take $[x]_\Pee\not= [y]_\Pee$ and $W \in \Pee$ such that $x\in W$ and $y\not\in W$. Fix sequences $ \{ U_n: n \in \omega\} $ and $ \{ V_n: n \in \omega\} $ witnessing $W \in  \Pee_{seq}$. Choose $k\in \omega$ such that $x\in U_k$ and $y\in V_k$. Hence  $ [x]_\Pee \subseteq U_k $  and $[y]_\Pee \subseteq V_k$. Therefore,  sets $q[U_k]$ and $ q[V_k]$  are disjoint neighbourhoods of  
$[x]_\Pee$ and $[y]_\Pee$, respectively.
 \end{proof} 
 \begin{lemma}\label{reg}  If a non-empty family of sets $\Pee \subseteq \Pee_{seq}$ is closed under finite intersections, then  $\Qee_\Pee$-topology  is   regular.
\end{lemma}
\begin{proof}
 We have $q[A] \cap q[B] = q[A\cap B]$ for each $A, B \in \Pee$. The family $\{ q[A]: A \in \Pee\}$ is a base of open sets for the $\Qee_\Pee$-topology. Fix $x\in W\in \Pee$ and sequences $ \{ U_n: n \in \omega\} \subseteq \Pee$ and $ \{ V_n: n \in \omega\} \subseteq \Pee$  witnessing $W\in \Pee_{seq}$. Take any $U_k \subseteq W$ such that  $[x]_\Pee \subseteq U_k\in \Pee$. We get   $q(x) \in q[U_k]\subseteq \cl q[U_k] 
 \subseteq q[X\setminus V_k] = X{/\Pee}\setminus q[V_k]\subseteq q[W]$, where $\cup \Pee=X$. 
  \end{proof}
To show which $\Qee_\Pee$-topologies are completely regular, we apply the Frink's theorem, compare \cite{fri} or \cite{eng} p. 72. 
   
\textbf{Theorem} [O. Frink (1964)]. \textit{ A $T_1$-space $X$ is completely regular if, and only if there exists  a base $\Bee$ satisfying}:\\
\indent (1) \textit{If  $x\in U\in \Bee$, then there exists  $V\in \Bee$ such that $x\not\in V$ and $U\cup V = X$};\\
\indent (2) \textit{ If   $U,V\in \Bee$ and 
 $U\cup V= X$, then there exists
disjoint sets $M,N \in \Bee$ such that $X\setminus U\subseteq M$ and   $X\setminus V\subseteq N$}.  \qed

  \begin{theorem}\label{t5} If $\Pee $ is a ring of subsets of $X$ and $\Pee \subseteq \Pee_{seq}$, then the $\Qee_{\Pee}$-topology is  completely regular.
\end{theorem}
\begin{proof} The $\Qee_\Pee$-topology is Hausdorff by Lemma \ref{h5}. Let $\Bee $ be the minimal family which contains
 $\{ q[V]: V \in \Pee \}$  and is closed under countable unions.  This family is a base for the $\Qee_\Pee $-topology, by Lemma \ref{tie}. We should show that $\Bee$ fulfills conditions $(1)$ and $(2)$ in  Frink's theorem. 

 Let $[x]_\Pee \in q[W]\in \Bee$. Fix sequences  $\{ U_k: k\in \omega\}$  and $\{ V_k: k\in \omega\}$  witnessing  $W\in  \Pee_{seq}$ and $k\in \omega$  such that $x\in X \setminus V_k \subseteq W$. We have $W \cup V_k = X$. Therefore $[x]_\Pee \not\in q[V_k]$ and $q[W] \cup q[V_k] = X/\Pee $. Thus $\Bee$ fulfills $(1)$.

Fix sets $\bigcup \{U_n: n\in \omega\}\in \Bee$ and $\bigcup \{V_n: n\in \omega\}\in \Bee$ such that   $$X/\Pee = \bigcup \{q[U_n]: n\in \omega\} \cup \bigcup \{q[V_n]: n\in \omega\},$$ where $U_n $ and  $ V_n$ belong to $ \Pee $.  Thus,     $U = \bigcup \{U_n: n\in \omega\}\in  \Pee_{seq}$ and  $V= \bigcup \{V_n: n\in \omega\} \in  \Pee_{seq}$ by Lemma \ref{4hau}. Next, fix sequences   $ \{ A_n: n \in \omega\} $,  $ \{ B_n: n \in \omega\} $, $  \{ C_n: n \in \omega\} $ and $ \{ D_n: n \in \omega\} $  witnessing $U \in \Pee_{seq}$ and $ V\in \Pee_{seq}$, respectively.  Therefore $$A_k\subseteq (X\setminus B_k)  \subseteq A_{k+1} \subseteq U   \mbox{ and } C_k\subseteq (X\setminus D_k)  \subseteq C_{k+1}\subseteq V,$$ for every $k\in\omega$. Put $N_n=A_n\cap D_n$ and $M_n=C_n\cap B_n$. Let $$M=\bigcup\{M_n:n\in\omega\} \mbox{ and } N=\bigcup\{N_n:n\in\omega\}. $$ Sets $q[M]$ and $q[N]$  fulfill $(2)$ in  Frink's theorem. 
Indeed, if $k\leq n$, then  $$A_k \cap D_k \cap C_n \cap B_n \subseteq A_n \cap B_n =\emptyset$$  and $$
A_n \cap D_n \cap C_k \cap B_k \subseteq C_n \cap D_n =\emptyset. $$ Consequently   $M_k\cap N_n=\emptyset$, for any $k,n\in \omega$. Hence   sets $q[M]$ and $q[N]$ are disjoint. 
Also, it is $q[V]\cup q[N]=X/\Pee $. Indeed, suppose that $x\not\in V$, then $x\in U$ and there is $k$ such that $x\in A_k$. Since $x\not\in V$, then $x\in D_k$ for all $k\in \omega$. We have $x\in A_k\cap D_k=N_k\subseteq N$. Therefore $[x]_\Pee  \in q[N]$. Similarly, one  gets  $q[U]\cup q[M]=X/\Pee $. Thus $\Bee$ fulfills $(2)$.
 \end{proof}
 
 If $ \Pee\subseteq \Pee_{seq}$ is finite, then $X/\Pee$ is discrete, being a finite  Hausdorff space. Whenever $\Pee\subseteq \Pee_{seq}$ is countable and closed under finite intersections, then  $X/\Pee$ is a  regular space with a countable base. Therefore,  $X/\Pee$ is  metrizable and separable.   
\section{Skeletal families and skeletal functions}
A continuous surjection  is called \textit{skeletal}  whenever for  any non-empty open sets $U\subseteq X$ the closure of  $f[U]$ has non-empty interior. If $X$ is a compact space and  $Y$  Hausdorff, then a continuous surjection $f:X\to Y$ is skeletal if, and only if $\int f[U] \not=\emptyset,$  for every non-empty and open $U\subseteq X$.   One can find  equivalent notions \textit{almost-open} or \textit{semi-open} in the literature,  see  \cite{arh} and \cite{her}. Following J. Mioduszewski and L. Rudolf \cite{mr} we call  such maps skeletal, compare \cite{s79} p. 413. In a fact, one can use the next proposition  as a definition for skeletal functions.    

\begin{proposition}\label{fer}
Let $f: X\to Y$ be a skeletal function. If  an open set $V \subseteq Y$ is dense, then the preimage $f^{-1} (V) \subseteq X$ is dense, too. 
\end{proposition} 
\begin{proof} Suppose that a non-empty open set $W\subseteq X$ is disjoint with $f^{-1} (V)$.  Then the image $\cl f[W]$ has non-empty interior and $\cl f[W]\cap V =\emptyset$, a contradiction.
\end{proof}

There are topological spaces with no skeletal map onto  a dense in itself metrizable space. For example, the remainder of the \v{C}ech-Stone
 compactification $\beta N$. Also, if $I$ is a compact segment of connected Souslin line and $X$ is metrizable, then each skeletal map $f:I\to X$ is constant. Indeed, let $Q$ be a countable and dense subset of   $f[I] \subseteq X$. Suppose a skeletal map $f:I\to X$ is  non constant. Then the preimage $f^{-1}(Q)$ is nowhere dense in $I$ as countable union of nowhere dense subset of a Souslin line. So, for each  open set $V\subseteq I\setminus f^{-1}(Q)$ there holds  $\int f[V] = \emptyset$, a contradiction. Regular Baire space X with a category measure $\mu $, for a definition of this space  see
\cite[pp. 86 - 91]{ox}, gives an another example of a space with no skeletal map onto  a dense in itself, separable and metrizable space.  
 In \cite{bs} A. B\l aszczyk
 and S. Shelah  are considered separable extremally disconnected spaces with no   skeletal map onto a dense in itself, separable and metrizable space. They formulated the result in terms of Boolean algebra:  \textit{There is a nowhere dense ultrafilter on $\omega$ if, and only if
there is a complete, atomless, $\sigma$-centered Boolean algebra which contains no regular, atomless,
countable subalgebra}. 

A family $\Pee$ of open subsets of a space $X$ is called a \textit{skeletal family}, whenever for every  non-empty open set $V\subseteq X$  there exists $W \in \Pee$ such that $U \subseteq W$ and $\emptyset\not= U \in \Pee$ implies $U \cap V \not= \emptyset $. The following proposition explains connection between skeletal maps and skeletal families.
\begin{proposition}\label{azs}
Let $f: X\to Y$ be a continuous function and let $\Bee$ be a $\pi$-base for $Y$. The family $\{f^{-1} (V): V \in \Bee\}$ is skeletal if, and only if $f$ is a skeletal map. 
\end{proposition} 
\begin{proof} Assume, that  $f$ is a skeletal map.  Fix a non-empty open set $V\subseteq X$. There exists $W\in \Bee$ such that $W\not=\emptyset$ and $W \subseteq \int
 \cl f[V]$. 
 Also, for any $U\in \Bee$ such that $\emptyset \not= U \subseteq W$ there holds $f^{-1} (U)\cap V \not= \emptyset$. Indeed,  if $f^{-1} (U)\cap V = \emptyset$, then  $U \cap \cl f[V] = \emptyset$, a contradiction. 
  Thus the family $\{f^{-1} (V): V \in \Bee\}$ is skeletal. 
 
Assume, that function $f:X\to Y$ is not skeletal. Then  there exists a non-empty open set $U\subseteq X$ 
such that $\int \cl f[U] = \emptyset$. Since $\Bee$ is a $\pi$-base for $Y$, then for each $W\in \Bee$ there exists $V\in \Bee$ such that $V\subseteq W$ and $V\cap f[U]=\emptyset$. The family $\{ f^{-1} (V): V \in \Bee \}$ is not skeletal.
\end{proof}
It is well know - compare  a comment following the definition of compact open-generated spaces in \cite{s81} - that   all limit projections are  open in any inverse  system  with open bonding maps. And conversely, if all limit projections of an inverse system are open,
then so are all bonding maps. Similar fact holds for skeletal maps.
\begin{proposition}\label{11}
If  $\{ X_\sigma, \pi^\sigma_\varrho,\Sigma\}$ is  a  inverse system such that  all bonding maps $\pi^\sigma_\varrho$ are skeletal and all projections  $\pi_\sigma$ are onto, then any projection $\pi_\sigma$ is skeletal.   
\end{proposition} 
\begin{proof} Fix $\sigma \in \Sigma$. Consider a  non-empty basic set $\pi^{-1}_\zeta (V)$ for the limit $\varprojlim \{ X_\sigma, \pi^\sigma_\varrho,\Sigma\}$. Take  $\tau \in \Sigma $  such that $\zeta \leq \tau $ and $  \sigma \leq \tau.$ 
We get $ \pi^{-1}_\zeta (V)= \pi^{-1}_\tau ( (\pi^{\tau}_\zeta)^{-1} (V)).$ Hence $$ \pi_\tau[\pi^{-1}_\zeta (V)] =  \pi_\tau[\pi^{-1}_\tau ( (\pi^{\tau}_\zeta)^{-1} (V))] = (\pi^{\tau}_\zeta)^{-1} (V),$$  the set $ \pi_\tau[\pi^{-1}_\zeta (V)]$ is open and non-empty. We have $$\pi_\sigma [\pi^{-1}_\zeta (V)]= \pi^\tau_\sigma [\pi_\tau[\pi^{-1}_\zeta (V)]]  ,$$ since $\pi^\tau_\sigma \circ \pi_\tau = \pi_\sigma$. The bonding map $\pi^\tau_\sigma$ is skeletal, hence  the closure of $\pi_\sigma [\pi^{-1}_\zeta (V)]$ has non-empty interior. 
 \end{proof}
 \section{The open-open game}    
    Players are playing at a topological space $X$ in the open-open game.   Player I  chooses  a non-empty open  subset $A_0\subseteq X$ at the beginning. Then Player II chooses a  non-empty open  subsets $B_0 \subseteq A_0$.   Player I  chooses 
a non-empty open  subset $A_n\subseteq X$  
at the $n$-th inning, and then  Player II chooses a non-empty open  subset $B_n \subseteq A_n$.
Player I wins, whenever  the union $B_0 \cup 
B_1 \cup  \ldots \subseteq X $ is  dense. One can assume that Player II wins for other cases. 
 The space $X$ is called  I-\textit{favorable} whenever Player I can be insured that he  wins no matter how Player II plays. In other words,  Player I has a winning strategy.   
A strategy for Player I could be  defined as a  function $$\sigma :\bigcup \{ \Tee^n: n\geq 0\} \to \Tee ,$$ 
where $\Tee$ is a family of non-empty and open subsets of $X$. 
Player I has a  
winning strategy, whenever  he knows how to define  $A_0=\sigma(\emptyset)$ and succeeding $A_{n+1}= \sigma( B_0,B_1 , \ldots, B_n)$ such that  for each game  $$ (\sigma (\emptyset), B_0, \sigma(B_0), B_1,\sigma(B_0, B_1),B_2, \ldots, B_n ,\sigma( B_0,B_1 , \ldots, B_n), B_{n+1},  \ldots )$$ the union $B_0 \cup B_1 \cup B_2  \cup  \ldots \subseteq X$ is  dense. For more details about the open-open game see P. Daniels, K. Kunen  and H.~Zhou \cite{dkz}. 
 
    Consider  a countable sequence $\sigma_0, \sigma_1, \ldots $ of strategies for Player I. For a  family $\Qee\subseteq\Tee$ let  $\Pee (\Qee ) $ be the minimal family   such that  $\Qee \subseteq \Pee(\Qee )\subseteq \Tee$, and if  $\{ B_0,B_1 , \ldots, B_n \} \subseteq \Pee(\Qee )$, then    $\sigma_k( B_0,B_1 , \ldots, B_n) \in \Pee(\Qee ),$  and $\sigma_k(\emptyset) \in \Pee(\Qee ),$ for all $\sigma_k$. We say that $\Pee(\Qee )$ is the \textit{closure of $\Qee$ under strategies }$\sigma_k$. In particular, if $\sigma$ is a  winning strategy and the closure of  $\Qee$ under $\sigma$ equals $\Qee$, then $\Qee $ is closed under a  winning strategy.

\begin{lemma}\label{l6} If   $\Pee $ is closed under a  winning strategy for Player I, then for any  open set $V\not=\emptyset$   there is $W \in \Pee$ such that whenever $U \in \Pee$ and $U \subseteq W$ then  $U \cap V \not= \emptyset $. 
\end{lemma}
\begin{proof}  Let $\sigma$ be a winning strategy for Player I. Consider an open set $V\not=\emptyset$. Suppose that for any $W \in \Pee$ there is  $U_W \in \Pee$ such that $U_W \subseteq W$ and  $U_W \cap V = \emptyset $. Then Player II wins any game whenever he always chooses  sets $U_W\in \Pee$, only.  In particular,  the game $$\sigma(\emptyset), U_{\sigma(\emptyset)}, \sigma( U_{\sigma(\emptyset)}), U_{\sigma( U_{\sigma(\emptyset)})}, \sigma (U_{\sigma(\emptyset)},U_{\sigma( U_{\sigma(\emptyset)})}), U_{\sigma (U_{\sigma(\emptyset)},U_{\sigma( U_{\sigma(\emptyset)})})},\ldots$$ would be winning for him, since all sets chosen by Player II: $$U_{\sigma(\emptyset)},  U_{\sigma( U_{\sigma(\emptyset)})},  U_{\sigma (U_{\sigma(\emptyset)},U_{\sigma( U_{\sigma(\emptyset)})})},\ldots ;$$ are disjoint with $V$,   a contradiction. 
 \end{proof}

\begin{theorem}\label{t7} If a ring  $\Pee$  of open subsets of $X$ is  closed   under a winning strategy and  $\Pee\subseteq  \Pee_{seq}$, then $X/{\Pee}$  is a completely regular space and the $\Qee_\Pee$-map $q:X\to X/{\Pee}$ is  skeletal.\end{theorem}
\begin{proof}
Take    a nonempty open  subset  $V \subseteq X$.    Since 
 $\Pee$ is closed under a winning strategy, there exists 
 $W \in \Pee$ such that  if $U \in \Pee$ and $U \subseteq W$, then  $U \cap V \not= \emptyset $, by Lemma \ref{l6}. This follows  
  $q [U] \cap q [V] \not= \emptyset $, for any basic set $q[U]$ such that $U\subseteq W$ and $U\in \Pee$. Therefore $q[W]\subseteq \cl q[V]$, since  $\{ q[U]: U \in \Pee\}$ is  a base for the $\Qee_\Pee$-topology.
The $\Qee_\Pee$-map $q:X\to X/{\Pee}$ is continuous 
by Lemma \ref{tie}. By Theorem \ref{t5}, the space $X/{\Pee}$  is  completely regular. 
\end{proof}
 
 Fix a $\pi$-base $\Qee$ for a space $X$. Following 
\cite{dkz}, compare \cite{kp},  any family $\Cee \subset [\Qee]^{\omega}$ is called \textit{a club filter} whenever:\\ \indent 
 The family $\Cee$ is closed under $\omega$-chains with respect to inclusion, i.e. if $\Pee_1 \subseteq \Pee_2 \subseteq  \ldots $ is an  $\omega$-chain which consists of elements of $\Cee$, then  $\Pee_1 \cup \Pee_2 \cup \ldots \in \Cee$;
 For any countable subfamily  $\Aaa\subseteq \Qee$, where $\Qee$ is the   $\pi-$base fixed above,  there  exists $\Pee \in \Cee$ such that   $ \Aaa\subseteq\Pee$; and

$(\See)$.  \textit{For any non-empty open set $V$ and each $\Pee \in \Cee$  there is $W \in \Pee$ such that if $U \in \Pee$ and $U \subseteq W$, then $U$ meets $V$, i.e. $U \cap V \not= \emptyset $}. 

In fact, the condition $(\See)$ gives reasons to look into I-favorable spaces with respect to skeletal families.  Any $\Pee $ closed under a  winning strategy for Player I fulfills $(\See)$, by Lemma \ref{l6}. There holds, see \cite{dkz} Theorem 1.6, compare \cite{kp} Lemmas 3 and 4: 
 \textit{A topological space has a club filter if, and only if it is I-favorable}.
 In the next part we modify a little the definition of  club filters. We introduce   $\Tee$-clubs, i.e.   club filters with some additional properties.

Suppose a completely regular space $X$ is I-favorable. Let  $\Tee$ be the family of all cozero subsets of $X$. 
  For each $W\in \Tee$ fix 
sequences   $ \{ U_n^W: n \in \omega\}$ and $ \{ V_n^W: n \in \omega\} $ witnessing $W\in \Tee_{seq}$.  First, for each $k$ choose $\sigma_k^*(\emptyset)\in \Tee$. Next, put $\sigma^*_{2n} (W) = U_n^W$ and $\sigma^*_{2n+1} (W) = V_n^W$, and  $\sigma_k^*(\See)=\sigma_k^*(\emptyset)$ for other $\See \in \bigcup \{ \Tee^n: n\geq 0\}$. Then,  a family  $\Pee \subseteq \Tee $ is
 closed  under strategies $\sigma_k^*$, whenever   $\Pee\subseteq \Pee_{seq} $. 
Also,  $\Pee$ is closed under finite unions, whenever it is closed under the strategy which  assigns the union $A_0 \cup A_1\cup \ldots \cup A_n$ to each sequence $ (A_0, A_1,\ldots , A_n) $. And also, $\Pee$ is closed under finite intersections, whenever it is closed under the strategy which assigns the intersection  $ A_0 \cap A_1\cap \ldots \cap A_n$ to each $   (A_0, A_1,\ldots , A_n)$.

 Consider a collection  $\Cee =\{ \Pee (\Qee ): \Qee \in [\Tee]^{\omega}\}$.   Assume that each $\Pee \in \Cee$ is countable and closed under a winning strategy for Player I and all strategies   $\sigma^*_k$,  and  closed under finite intersections and finite unions. Then, the family $\Cee$ is called  $\Tee$-\textit{club}. 
By the definitions, any  $\Tee$-\textit{club} $\Cee$   is closed under $\omega$-chains with respect to the inclusion. Each $\Pee\in \Cee$ is  a countable ring of sets and $\Pee\subseteq \Pee_{seq}$ and it is closed under a winning strategy for Player I. By Theorem \ref{t7},  the $\Qee_\Pee$-map $q:X\to X/\Pee$ is   skeletal and  onto a metrizable separable space, for every  $\Pee \in \Cee$.    
 
 Thus, we are ready  to build an inverse system with skeletal bonding maps onto metrizable separable spaces. 
Any $\Tee$-club $\Cee$ is  directed by the inclusion. For each $\Pee\in\Cee$ it is assigned the space $X/\Pee$ and the skeletal function  $q_\Pee:X\to X/\Pee$. If $\Pee, \Raa\in\Cee$ and $\Pee \subseteq \Raa$, then put $q^\Raa_\Pee ( [x]_\Raa) = [x]_\Pee.$ Thus, we have defined  the inverse system $\{ X/\Raa, q^\Raa_\Pee,\Cee\}$. Spaces $X/\Raa$ are metrizable and separable, bonding maps $q^\Raa_\Pee$ are skeletal and the directed set $\Cee$ is $\sigma$-complete. 
\begin{theorem}\label{9}
Let $X$ be a I-favorable completely regular  space. If $\Cee$ is  a $\Tee$-club, then the limit $ Y= \varprojlim \{ X/\Raa, q^\Raa_\Pee,\Cee\}$ contains a dense subspace which is homeomorphic to $X$. 
\end{theorem}
\begin{proof} For any $\Pee\in\Cee$, put $f(x)_\Pee=q_\Pee(x)$.
 We have defined the function  $f:X\to Y$ such that $f(x) = \{f(x)_\Pee\}$. If  $\Raa, \Pee \in \Cee$ and  $\Pee \subseteq \Raa$, then  $q^\Raa_\Pee (f(x)_\Raa)=f(x)_\Pee$. Thus  $f(x)$ is a thread, i.e. $f(x)\in Y$. 
 
The function $f$ is continuous. Indeed,  let $\pi_\Pee$ be the projection of $Y$ to $X/\Pee$. By \cite{eng} Proposition 2.5.5, the family  $\{\pi^{-1}_\Pee(q_\Pee[U]):U\in\Pee\in\Cee \}$ is a base for $Y$.  Also,   
$$f^{-1}(\pi^{-1}_\Pee(q_\Pee[U]))=q^{-1}_\Pee(q_\Pee[U])=U$$ holds for any $U\in\Pee\in\Cee  $. 

Verify that $f$ is injection.
Let $x,y\in X$ and $x\not =y$. Take $\Pee\in\Cee$ such that   $x\in U$ and 
$y\in V$ for some disjoint sets $U,V\in\Pee$. Sets $q_\Pee[U]$ and $q_\Pee[V]$ are disjoint, hence $\pi^{-1}_\Pee(q_\Pee[U])$ and $\pi^{-1}_\Pee(q_\Pee[V])$ are disjoint neighbourhoods of  $f(x)$ and 
$f(y)$, respectively. 

There  holds $\;f[U]=f[X]\cap\pi^{-1}_\Pee(q_\Pee[U])$, whenever $U\in\Pee\in\Cee$. Indeed, $f[U]\subseteq \pi^{-1}_\Pee(q_\Pee[U])$ implies   $f[U]\subseteq f[X]\cap\pi^{-1}_\Pee(q_\Pee[U]).$ Suppose,  
there exists $y\in \pi^{-1}_\Pee(q_\Pee[U])\cap f[X]$ such that $y\not\in  f[U])$. Take $x\in X$ such that $f(x)=y$ 
and $x\not\in U$.  
 We get  $\pi_\Pee(f(x))=q_\Pee(x)\not\in q_\Pee[U]$, but this follows  $f(x)\not\in \pi^{-1}_\Pee(q_\Pee[U])$, a 
contradiction.

 Thus, $f$ is open, since $\Tee = \bigcup\Cee$ is a base for $X$. But $f[X]\subseteq Y$ is dense, since the family  $\{\pi^{-1}_\Pee(q_\Pee[U]):U\in\Pee\in\Cee \}$ is a base for $Y$. 
\end{proof}
\section{Reconstruction of  I-favorable spaces}

Now, we are ready to prove the announce analog of Shchepin's openly generated spaces. 

\begin{theorem}\label{10}
If $X$ is a I-favorable compact  space, then  $$ X = \varprojlim \{ X_\sigma, \pi^\sigma_\varrho,\Sigma\},$$ where   $\{ X_\sigma, \pi^\sigma_\varrho,\Sigma\}$ is a  $\sigma$-complete inverse system, all spaces $X_\sigma$ are compact and metrizable, and all bonding maps $\pi^\sigma_\varrho$ are skeletal and onto.  
\end{theorem} 
\begin{proof}
Let  $\Cee$ be  a $\Tee$-club. Put  $$ \{ X_\sigma, \pi^\sigma_\varrho,\Sigma\}= \{ X/\Raa, q^\Raa_\Pee,\Cee\}.$$   Each space $X_\sigma=X/\Raa$ has countable base, by the definition of $\Tee$-club. Also, each  $\Qee_\Raa$-map $q_\Raa: X \to X/\Raa$ is continuous, by Lemma \ref{tie}. Hence, any  space $X_\sigma$ is  compact and metrizable,   by Lemma \ref{reg}.  Each $\Qee_\Raa$-map  $q_\Raa:X\to X_\sigma$ is skeletal, by Theorem \ref{t7}. Thus, all bonding maps $\pi^\sigma_\varrho$ are skeletal, too.  The space $X$  is homeomorphic to a dense subspace of   $\varprojlim \{X_\sigma, \pi^\sigma_\varrho,\Sigma\}$, by Theorem \ref{9}.   We get  $ X= \varprojlim \{ X_\sigma, \pi^\sigma_\varrho,\Sigma\},$ since $X$ is compact.

The inverse system $\{X_\sigma, \pi^\sigma_\varrho,\Sigma\}$ is $\sigma$-complete. Indeed, suppose that  $\Pee_0 \subseteq \Pee_1 \subseteq \ldots$ and all $\Pee_n\in \Cee$. Let $\Pee= \bigcup \{\Pee_n : n\in \omega\} \in \Cee$. Put  $$ (h([x]_\Pee))_{\Pee_n} = q^{\Pee}_{\Pee_n}([x]_\Pee) = [x]_{\Pee_n}.$$ Since maps $q^{\Pee}_{\Pee_n}$ are continuous, we have defined a continuous function  $h:X/\Pee \to \varprojlim \{ X/\Pee_n, q^{\Pee_{n+1}}_{\Pee_n}\}.$  Whenever $\{ [x_n]_{\Pee_n} \}$ is a thread in the inverse system $\{X/\Pee_n,q^{\Pee_{n+1}}_{\Pee_n} \}$, then there exists   $x\in \bigcap\{ [x_n]_{\Pee_n}: n \in \omega\} $, since sets $[x_n]_{\Pee_n}$ consists of a centered family  of nonempty closed sets in a compact space $X$. Thus $h^{-1}(\{ [x_n]_{\Pee_n} \})= [x]_\Pee \in X/\Pee$, hence $h$ is a bijection.
\end{proof}

To obtain the converse of Theorem \ref{10} one should consider an inverse system of compact metrizable spaces with all bonding maps skeletal. Such assumptions are  unnecessary. So, we assume that spaces $X_\sigma$ have countable $\pi$-bases, only. 

\begin{theorem}\label{111}
Let  $\{ X_\sigma, \pi^\sigma_\varrho,\Sigma\}$ be a $\sigma$-complete inverse system such that  all bonding maps $\pi^\sigma_\varrho$ are skeletal and all projections $\pi_\sigma$ are onto. If all spaces $X_\sigma$ have countable $\pi$-base, then the limit $\varprojlim \{ X_\sigma, \pi^\sigma_\varrho,\Sigma\}$ is I-favorable.
\end{theorem}
\begin{proof} Let $\leq$ denotes the relation  which directs $\Sigma$.
 Describe the following strategy for a match playing at the limit $X = \varprojlim \{ X_\sigma, \pi^\sigma_\varrho,\Sigma\}$. Assume that Players play with basic sets of the form $\pi_\sigma^{-1}(V)$, where $V$ is  non-empty and open in $X_\sigma$ and  $\sigma\in \Sigma$.

Player I  chooses an open non-empty set $A_0\subseteq X$ at the beginning. Let $\Bee_0 =\{B_0\}$ be a respond of Player II. Take $\sigma_0\in\Sigma$ such that $B_0=\pi_{\sigma_0}^{-1}(V^0_0)\subseteq A_0$. Fix a countable $\pi$-base $\{V^0_0,V^0_1,\ldots\}$ for $X_{\sigma_0}$. 

Assume, that we have just settled indexes $\sigma_0\leq\sigma_1\leq\ldots \leq\sigma_n$ and $\pi$-bases $\{V^k_0,V^k_1,\ldots\}$ for $X_{\sigma_k}$, where $0\leqslant k\leqslant n$.   Additionally assume, that for any $V^k_m$ there exists $V^{k+1}_j$ such that 
$\pi^{-1}_{\sigma_{k+1}}(V^{k+1}_j) = \pi^{-1}_{\sigma_{k}}(V^{k}_m)$. Now, Player I plays each set from $$\Aaa_{n+1} =\{ \pi^{-1}_{\sigma_{k}}(V^{k}_m): k\leqslant n \mbox{ and } m\leqslant n\}$$ one after the other. Let $\Bee_{n+1}$ denote the family of all responds of Player II, for innings from $\Aaa_{n+1}$. Choose $\sigma_{n+1} \geq \sigma_n$ and  a countable $\pi$-base $\{V^{n+1}_0,V^{n+1}_1,\ldots\}$  for $X_{\sigma_{n+1}}$ which contains the family $$\{ (\pi^{\sigma_{n+1}}_{\sigma_{k}})^{-1}(V^{k}_m): k\leqslant n \mbox{ and } m \in \omega\}$$ and such that for any $V \in \Bee_{n+1}$ there exists $V^{n+1}_j$ such that 
$\pi^{-1}_{\sigma_{n+1}}(V^{k+1}_j) = V$.

  Let $\sigma = \sup \{ \sigma_n: n \in \omega \}\in \Sigma .$ Any set $\pi_{\sigma_n} [\bigcup\{\bigcup \Bee_n: n \in \omega\}]$ is dense in $X_{\sigma_n}$, since it intersects any $\pi$-basic set $V^n_j \subseteq X_{\sigma_n}$. The inverse system is $\sigma$-complete, hence the set $\pi_\sigma[\bigcup\{\bigcup \Bee_n: n \in \omega\}]$ is dense in $X_\sigma$.  The projection $\pi_\sigma$ is skeletal by Proposition \ref{11}. So, the set $\bigcup\{\bigcup \Bee_n: n \in \omega\}$ is  dense in $X$ by Proposition \ref{fer}. 
 \end{proof}
A continuous and open map is skeletal, hence every compact openly generated  space is I-favorable. 
\begin{corollary}\label{21}
Any compact openly generated  space is  I-favorable. \hfill $\Box$
\end{corollary} 

The converse is not true.  For instance,  the \v{C}ech-Stone compactification $\beta N$ of positive integers with the discrete topology  is I-favorable and extremally disconnected. But $\beta N$ is not openly generated, since a compact extremally disconnected and openly generated space has to be discrete, see Theorem 11 in \cite{s76}.

\noindent \textbf{Acknowledgement}

The authors wish to thank to  referees for  their careful reading of a first version of this paper and for comments
that have been very useful to improve the final form of the proofs of some results.

\end{document}